\documentclass[preprint,12pt,3p]{elsarticle}
\usepackage{amsmath}
\usepackage{amssymb}
\usepackage{mathtools}
\usepackage{physics}
\usepackage{graphicx}
\usepackage[colorlinks,citecolor=blue]{hyperref}

\journal{Journal of Computational Physics}

\begin{document}

\begin{frontmatter}

\title{Spectral shock detection for dynamically developing discontinuities}

\author[label1,label2,label3]{Joanna Piotrowska\corref{cor1}}
\address[label1]{Center for Nonlinear Studies, Los Alamos National Laboratory, Los Alamos, NM 87545, USA}
\address[label2]{Kavli Institute for Cosmology, University of Cambridge, Madingley Road, Cambridge, CB3 0HA, UK}
\address[label3]{Cavendish Laboratory, University of Cambridge, 9 JJ Thompson Avenue, Cambridge, CB3 0HE, UK}

\cortext[cor1]{I am the corresponding author}
\ead{jmp218@am.ac.uk}

\author[label1,label4,label5]{Jonah M. Miller}
\address[label4]{CCS-2, Los Alamos National Laboratory, Los Alamos, NM 87545, US}
\address[label5]{Center for Theoretical Astrophysics, Los Alamos National Laboratory, Los Alamos, NM 87545, US}
\ead{jonahm@lanl.gov}

\begin{abstract}
Pseudospectral schemes are a class of numerical methods 
capable of
solving smooth problems with high accuracy thanks to their
exponential convergence to the true solution.
When applied to discontinuous problems,
such as fluid shocks and material interfaces, due to the Gibbs
phenomenon,
pseudospectral solutions lose their superb convergence and suffer
from spurious oscillations across the entire computational domain.
Luckily, there exist theoretical
remedies for these issues which have been 
successfully tested in practice
for cases of well defined discontinuities.
We focus on one piece of this procedure---detecting a discontinuity
in spectral data. We show that realistic applications 
require treatment of discontinuities	
dynamically developing in time and that it poses 
challenges associated
with shock detection. More precisely, smoothly steepening
gradients in the solution spawn spurious oscillations due to
insufficient resolution, causing premature shock identification
and information loss.
We improve existing spectral shock detection techniques
to allow us to automatically detect true discontinuities
and identify cases for which post-processing is required to
suppress spurious oscillations resulting from the loss of resolution.
We then apply these techniques to solve an inviscid Burgers'
equation in 1D, demonstrating that our method correctly treats
genuine shocks caused by wave breaking and removes oscillations
caused by numerical constraints.
\end{abstract}


\end{frontmatter}


\section{Introduction}
\label{sec:introduction}

Pseudospectral methods are a powerful numerical approach, offering
exponential convergence to the true solution across the entire
domain for smooth functions (see \cite{Boyd2013} for a~review). 
Thanks to their exquisite accuracy and efficiency they have been
successfully applied in, for example, large scale relativistic 
astrophysics simulations (e.g. \cite{SPeC, Boyle2008, Scheel2009}). 
A~wealth of astrophysical phenomena, however, 
involve discontinuities such as fluid shocks in supernova explosions
or galactic spiral arms. For this class of problem, pseudospectral
schemes perform poorly due to sublinear convergence properties and 
spurious oscillations induced in the solution known as the 
\textit{Gibbs phenomenon} \cite{Wilbraham1848, Gibbs1898,
Michelson1898}. 

A substantial body of research has been conducted in the attempt to
alleviate the Gibbs phenomenon and recover the exponential 
convergence of pseudospectral schemes. Some proposed solutions
are reprojection onto combinations of different basis sets \cite{Gottlieb1992, Gottlieb2011}, 
filtering \cite{Vandeven1991}, artificial viscosity 
\cite{Tadmor1990}, modelling and removal of oscillations 
\cite{Krylov2006, Lipman2010} and mollification \cite{Gottlieb1985}. 
In \cite{Piotrowska2019} we explored and numerically tested the last
approach, using an optimal mollifier as prescribed by 
\cite{Tanner2006}. We showed successful recovery
of convergence properties away from discontinuities and offered
improvements to preserve the discontinuous character of the solution.
We also solved a toy problem advecting a discontinuous function, 
demonstrating greater robustness of mollification as compared to the
more popular Gegenbauer reconstruction methods. 

In this work we take the next step towards applying pseudospectral
methods to solve physically relevant shock problems. This time
we allow the discontinuities to develop dynamically, instead
of introducing them by construction and apply
previously explored methods to detect the created shocks. Robust 
identification of evolving discontinuities is crucial in 
pseudospectral shock capturing, regardless of the choice of
Gibbs phenomenon removal method.
We discuss challenges
associated with the detection of discontinuity onset and capturing
of steep gradients at constant spatial resolution. We
also present our heuristic solution to those issues, offering 
practical improvements to shock detection methods in case of shocks
dynamically developing during the simulation. 
Finally, we solve
a toy problem of the inviscid Burgers' equation to demonstrate the
performance of our methods when applied to shocks developing smoothly
over time.

\vspace{1cm}
\section{Background}
\label{sec:background}

In this section we review existing methods which provide basis
for our improvements presented further in Section~\ref{sec:our-improvements}.

\subsection{Overview of pseudospectral methods used in this work}
\label{sec:overview}

The pseudospectral decomposition is a means of discretizing 
a function
of interest with a~finite sum of $N$ orthogonal basis functions 
$T(x)_n$:
\begin{equation}
	u(x) \approx S_N[u](x)=\sum_{n=0}^N \hat{u}_n T_n(x)\, ,
\end{equation}
where $S_N[u]$ is the partial sum of $u$ and the basis functions
are denoted $T_n(x)$. The expansion coefficients $\hat{u}_n$ are chosen such that
\begin{equation}
  \lim_{N\to\infty} S_N[u](x) = u(x)
\end{equation}
and are given by:
\begin{equation}
	\hat{u}_n = \frac{(u,T_n)}{(T_n, T_n)}
\end{equation}
where
\begin{equation}
	(u, T_n) = \int_\Omega u(x) T_n(x) w(x)dx	 
\end{equation}
is the inner product between $u$ and $T_n$ on a domain $\Omega$. 
The integral  can be approximated via Gauss quadrature in which
case the $x_i$ are the quadrature collocation points and $w_i$ their
associated weights, such that:
\begin{equation}
\lim_{N \to \infty}  \sum_{i=0}^{N} u(x_i) T_n(x_i) w_i = \int_\Omega u(x) T_n(x) w(x)dx	 \, .
\end{equation}

Thus, one can conveniently transform the
collocation-point values $u(x_i) \equiv u_i$ into expansion 
coefficients $\hat{u}_n$ through\footnotemark:
\begin{equation}
	\hat{u}_n = \mathcal{V}_{ni} u_i\, ,
	\quad 
	\mathcal{V}_{ni} \equiv \frac{T_n(x_i)w_i}{(T_n, T_n)}
\end{equation}
where $\mathcal{V}_{ni}$ is the inverse of the 
\textit{Vandermonde matrix}. One then can also approximate the
partial sum of the derivative of $u$: $S_N[\partial_x u]$ with the
following differentiation matrices:
\begin{equation}
	\mathcal{M}_{mn}=\frac{(\partial_x T_m, T_n)}{(T_n, T_n)}\,,
	\quad
	\mathcal{N}_{ji}=\mathcal{V}_{jn}^{-1} \mathcal{M}_{nm} 
	\mathcal{V}_{mi}\,,
\end{equation}
such that differentiation conveniently reduces to the following 
matrix operations:
\begin{equation}
	\lim_{N \to \infty} \dv{x}\hat{u}_m = \mathcal{M}_{mn}\hat{u}_n\,,
	\quad
	\lim_{N \to \infty} \dv{x}u_j = \mathcal{N}_{ji} u_i\, .
\end{equation}

\footnotetext{where $\sum_i a_i b_i = a_i b_i$ in the Einstein
notation}

In our work we use a basis of Chebyshev polynomials of the 
first kind\footnote{other choices are possible, depending on the
problem of interest e.g. a Legendre polynomial or Fourier basis} 
defined on the domain $\Omega = [-1,1]$ and 
the \textit{Chebyshev-Gauss-Lobatto} quadrature.
The associated grid of collocation points allows us to probe the
solution values at the boundaries with $x_0 = -1$ and $x_N = 1$.
For a detailed
description of Gaussian quadrature and transformations between
the spectral and collocation-point bases we refer the interested 
reader to \cite{Grandclement2006} and \cite{Boyd2013}.

\vspace{1cm}
\subsection{Stabilization techniques for non-linear PDEs}
\label{sec:stabilizing}

In this research we focus on one-dimensional (1D) problems which 
smoothly develop
discontinuities over time, such as the non-linear \textit{1D
inviscid Burgers' equation}:
\begin{equation}
	\partial_t u + u\, \partial_x u = 0.
	\label{eq:inviscid-burgers}
\end{equation} 

In our investigation we assume smooth initial conditions, 
choosing a Gaussian $g(x)$ centered on $x=0$:
\begin{equation}
	g(x) = \exp(\frac{-(x-x_0)^2}{2\sigma^2})\, ,
	\label{eq:gaussian}
\end{equation}
where $x_0 = 0$ and $\sigma=0.15$.

Once we discretize the initial conditions with a pseudospectral
projection onto the Chebyshev polynomial basis, our hyperbolic partial
differential equation (PDE) reduces to a set of $N+1$ ordinary
differential equations (ODEs) at the $N+1$ collocation points. Thus,
we solve Eq.~(\ref{eq:inviscid-burgers}) using the method of lines,
keeping time $t$ as a continuous variable and assuming periodic
boundary conditions such that at each time step $u_0=u(-1)$ is set to
$u_N=u(1)$. Note that our initial condition is non-smooth at the
boundaries, since $g(x)$ is non-periodic. Our periodic boundary
conditions serve as a simplistic proxy for a multi-domain spectral
method, with the left and right boundaries serving as inflow and
outflow into other spectral domains respectively. This allows us to
run our simulation for arbitrary long time, despite the finite-sized
domain.

In general, any linearly stable numerical scheme 
does not guarantee stability for non-linear problems.
In a pseudospectral method
nonlinear terms can produce artificially large short-wavelength
coefficients which need to be subdued to avoid a numerical
catastrophe.
In this work we choose a modification to the spectral viscosity 
technique \cite{Tadmor1990}, namely the right hand side filtering
introduced by \cite{Miller2016} in Eq.~(73-75). In this framework
filtering is applied to the derivative operator modes via the
following diagonal matrix:
\begin{equation}
	\mathcal{F}_{mn}:= \qty(\frac{n}{N})^{2s} \delta_{mn}\, ,
\end{equation}
where $s$ is the so-called order of dissipation and we chose $s=2$
in the evolution of our system which proved optimal in our numerical
tests. The temporal derivative at each grid point then takes the
form:
\begin{equation}
	\partial_t u_i = u_i\, \mathcal{N}_{ij}u_j
	- c\, \mathcal{V}_{im}^{-1} \mathcal{F}_{mn} 
	\mathcal{V}_{nj} u_j\, ,
\end{equation}
where the constant $c$ depends on the particular problem and 
$c \sim 0.01$ performed best when tested with the inviscid Burgers' 
equation.  

\vspace{1cm}
\subsection{Removing the Gibbs phenomenon}
\label{sec:Gibbs}

The main advantage of pseudospectral methods is their
\textit{global exponential convergence} for smooth functions. 
More specifically, the exponential decrease of the 
interpolation error $u - S_N[u]$
integrated across the computational domain
as a function of expansion order $N$:
\begin{equation}
	\norm{u - S_N[u]}_2 \leq \frac{\alpha}{N^p} 
	\sum_{k=0}^p \norm{u^{(k)}}_2\, ,
\end{equation}
where $0 \leq p \leq p_{max}$, $p_{max}$ is the highest order of
derivative well-defined for $u$, $\alpha$ is a constant factor and
$\norm{f(x)}_2 \equiv \int_\Omega \abs{f(x)}^2 dx$ denotes the $L_2$
norm. In the smooth case, arbitrarily high-order derivatives are
well-defined and the inequality holds for all $p > 0$
\cite{Grandclement2009}.

\begin{figure}
	\includegraphics[width=\textwidth]{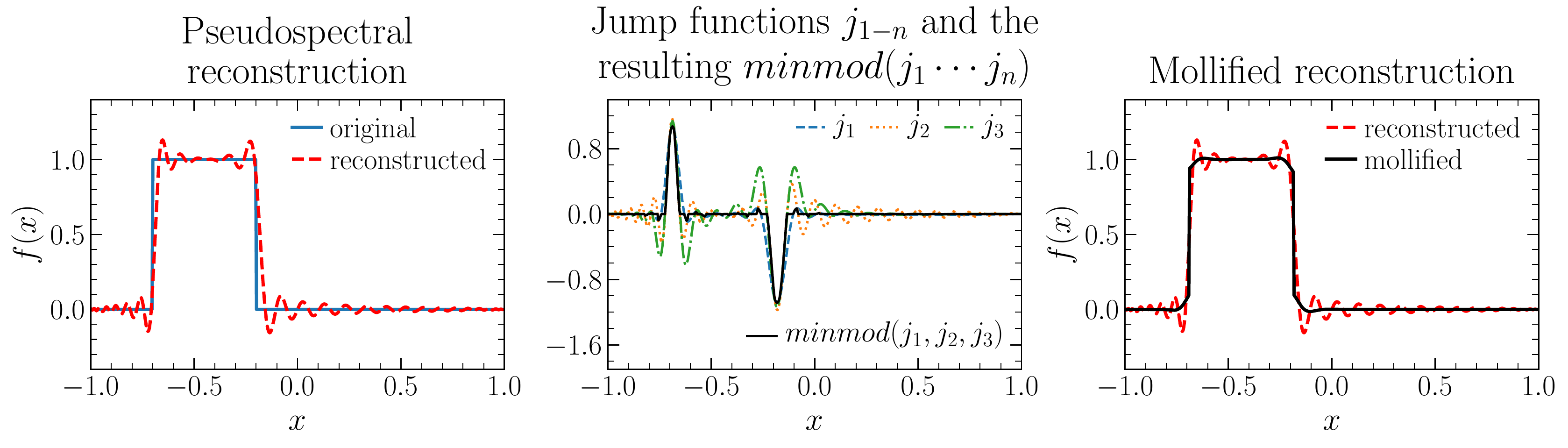}
	\caption{Left panel: demonstration of the Gibbs oscillations
	affecting the pseudospectral reconstruction (red dashed line)
	of a discontinuous `tophat' function (Eq.~(\ref{eq:tophat}), 
	blue solid line).
	Middle panel: exemplary jump functions and the resulting
	$minmod$ function from \cite{Gelb2008}. The extrema indicate
	locations of the discontinuities in the `tophat'.
	Right panel: mollified reconstruction (black solid line) 
	post-processed using one-sided mollifiers from 
	\cite{Piotrowska2019}. The reconstruction is free from
	the oscillations while preserving the sharp character of
	its discontinuities. This illustration was made for 
	expansion order $N=60$.}
	\label{fig:demo}
\end{figure}

When applied to discontinuous
problems, pseudospectral interpolants lose their convergence
properties with the
$L_2$ norm exhibiting at best $\mathcal{O}\qty(N^{-\flatfrac{1}{2}})$
convergence \cite{Thomson2008}.
Furthermore, the pseudospectral representations of the solution
suffer from spurious oscillations. 
The 
oscillation amplitude decreases away from the discontinuities and
their frequency varies depending on the position within the domain.
The left panel of Fig.~\ref{fig:demo} 
illustrates this phenomenon by overlaying a Gibbs-affected $N=60$ 
reconstruction on the input discontinuous piecewise linear 
`tophat' function:
\begin{equation}
	h(x) = 
	\begin{dcases}
		1 & \text{if } x_1 < x < x_2 \\
		0 & \text{otherwise}
	\end{dcases}\, , 
	\label{eq:tophat}
\end{equation}
where $x_1 = -0.7$ and $x_2=-0.2$.

In order to mitigate these issues and recover a Gibbs-free solution
one needs to post-process the affected spectral representation. 
This procedure consist of two main parts: \textit{edge detection}
and \textit{Gibbs oscillations removal}.
Edge detection is a step necessary for many Gibbs phenomenon
removal techniques,
allowing construction of adaptive smoothing kernels 
in the \textit{mollification} approach or the correct domain 
decomposition in the case of the Gegenbauer reconstruction 
\cite{Gottlieb2011a}.

Following our work in \cite{Piotrowska2019}, we choose to 
locate the edges using the $minmod$ technique introduced
by \cite{Gelb2008}. It relies on applying a range of different
concentration factors
$\mu(\frac{k}{N})$ \cite{Gelb1999} to the partial
sum of the function derivative $\tilde{S}_N u ^\prime (x)$, which 
accelerate the rate at which $\tilde{S}_N u ^\prime (x)$ 
converges onto the \textit{jump function},
\begin{equation}
  [u](x) = u(x^+)-u(x^-).
\end{equation}
In this prescription, each choice of an admissible 
$\mu(\frac{k}{N})$ generates a jump function approximation
$j_\mu(x)$ given by the following partial sum:
\begin{equation}
	j_\mu (x) = {\frac{\pi \sqrt{1-x^2}}{N} 
	\sum^{N}_{k=1} \mu\qty({\frac{k}{N}}) 
	\hat{u}_k T^\prime_k(x)}\, .
	\label{eq:jump}
\end{equation}


The jump function 
approximation then converges onto the true jump function 
at the following rate \cite{Gelb2008}:
\begin{equation} 
	\abs{j_\mu (x) - [u](x)}\leq 
	Const \cdot \frac{\log{N}}{N}
\end{equation}
and each $j_\mu(x)$ exhibits a different oscillation pattern
around the discontinuity associated with a given 
$\mu(\frac{k}{N})$ \cite{Gelb2000}. The $minmod$ prescription
then takes advantage of the differences in the oscillatory patters, 
combining a range of jump functions such that the oscillations
are removed \cite{Gelb2008}:
\begin{equation}
    \label{eq:minmod}
    minmod(j_1,\ldots,j_n)(x):=
    \begin{cases}
        \min(j_1(x),\ldots,j_n(x)), & \text{if $j_1(x),\ldots,j_n(x)>0$,} \\
        \max(j_1(x),\ldots,j_n(x)), & \text{if $j_1(x),\ldots,j_n(x)<0$,} \\
        0, & \text{otherwise.}
    \end{cases}
\end{equation}

The middle panel of Fig.~\ref{fig:demo} demonstrates three jump 
function approximations
$j_\mu(x)$ from trigonometric, polynomial and exponential 
$\mu(\frac{k}{N})$ forms as well as a resulting $minmod(x)$
for $N=60$.
The location of discontinuities is indicated by the positions of
the peaks in the function represented by the black solid curve.
In our work we construct the $minmod$ using $\mu(\frac{k}{N})$ 
defined in Eq.~(8-10) in \cite{Gelb2008} and their filtered versions
constructed using Lanczos \cite{Lanczos1966} filters of 
orders 0, 1, 2 and 3. 

Once we identify the discontinuities and their positions we
then remove the oscillations through mollification
with adaptive one-sided mollifiers \cite{Tanner2006, 
Piotrowska2019}. The
right panel of Fig.~\ref{fig:demo} demonstrates the result of post
processing by overplotting the mollified result (black solid) on 
the Gibbs-affected reconstruction (red dashed). The post-processed
solution is successfully freed from spurious oscillations,
simultaneously preserving its discontinuous character.

Finally, it is important to mention, that even in the ideal
cases of well-defined discontinuities the $minmod$ function is 
polluted by residual oscillations of small amplitude which can
be misidentified as jumps in an automated search for extrema. 
Hence one always has to make a choice of a peak height qualifying
as a true extremum. 
As a result, one then sets the range of discontinuity heights
captured by the edge detection framework, which depends on the
resolution and the problem at hand.

\vspace{0.8cm}
\section{New shock detection challenges}
\label{sec:detection-challenges}

The edge detection techniques illustrated in Fig.~\ref{fig:demo}
are an implementation of prescriptions developed by A. Gelb and
collaborators \cite{Gelb2000, Gelb2008} for stationary 
discontinuities.
Realistic numerical applications, however, involve dynamical 
development of discontinuities in time and thus require the
pseudospectral framework to capture transitions between the 
smooth and
discontinuous (shocked) regimes in an automated fashion. Thus, one 
requires a~robust means of automatically recognizing discontinuous 
functions to deploy post-processing techniques at appropriate 
time steps.

\begin{figure}
	\includegraphics[width=\textwidth]{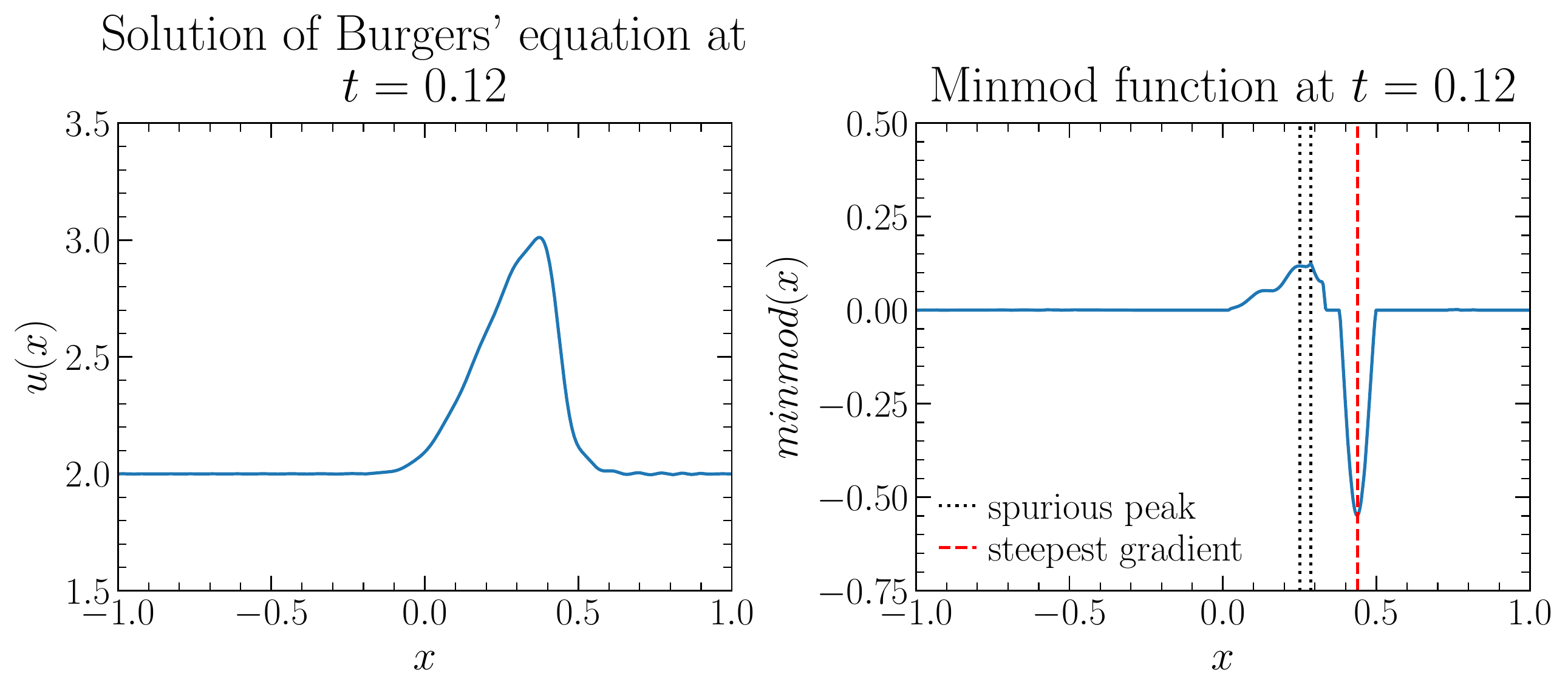}
	\caption{Issues arising in discontinuities developing
	dynamically from smooth solutions. Left panel: mild
	oscillations in the solution present in the vicinity 
	of steep gradients. The pseudospectral basis fails
	to capture overly steep gradients due to its
	resolution limit. 
	Right panel: corresponding $minmod$ function polluted
	by spurious peaks. These are misidentified as shock
	locations in an automated search.}
	\label{fig:challenges}
\end{figure}

Solutions in which shocks develop smoothly over time experience
additional problems caused by insufficient resolution of the 
truncated pseudospectral expansion basis. An example of such is
the inviscid Burgers' equation (eq.~\ref{eq:inviscid-burgers}) 
where the spatial gradient in the function
steepens continuously until the point of wave breaking. In such 
cases the solution is populated by additional oscillations around 
steep gradients which result from the lack of sufficient 
resolution to capture significant variation at very small scales. 
The left panel 
of Fig.~\ref{fig:challenges} demonstrates this issue by showing
the affected solution $u(x)$ of the 1D inviscid Burgers' 
equation \ref{eq:inviscid-burgers} with smooth Gaussian initial
conditions (Eq.~\ref{eq:gaussian} with $x_0=0$ and $\sigma=0.15$),
subject to periodic boundary conditions. 
The panel shows a skewed 
Gaussian at time step $t=0.12$ which has a
steep gradient at $x\approx0.45$. The solution in the vicinity of 
the
steep slope in $u(x)$ is populated by spurious oscillations of small
amplitude which introduce undesired, incorrect, information. These
oscillations are an artefact of insufficient resolution and their
amplitude converges away with increasing expansion order $N$.


What is more, the resolution limit also demonstrates itself in 
shock detection, altering the $minmod$ function such that it
exhibits different convergence properties and features additional 
peaks. The right panel of Fig.~\ref{fig:challenges} highlights the
latter by showing the irregular shape of the $minmod$ to the left
of the clear peak corresponding to a steep gradient in $u(x)$. 
The black dotted lines indicate the location of two spurious peaks 
identified in an automated search, which would be wrongly 
interpreted as positions of two discontinuities within the 
domain. This flavor of deformations present in the $minmod$ due to
insufficient resolution causes the framework to detect 
discontinuities before they have fully formed and wrongly
identify locations of shocks which are not present in the solution.
As in the case of spurious oscillations in the solution, the 
$minmod$ artefacts depend on the expansion order $N$ and disappear
at sufficiently high values of $N$.

In summary, dynamically developing shocks interfere with the 
automated
detection of discontinuities through the stencil's inability to
capture steep gradients. This results in premature shock
identification and gives rise to localized Gibbs-like oscillations in
the vicinity of unresolved gradients.

\vspace{0.8cm}
\section{Shock detection improvements}
\label{sec:our-improvements}

In order to address the described issues we first check whether
at a given time step the solution is continuous within our resolution
limit. If that is not the case, we then determine whether the 
solution requires treatment appropriate for a genuine discontinuity
or whether the resolution is not sufficient to capture steep 
gradients in the function.

\subsection{Identification of resolved continuous cases}
\label{sec:smooth-resolved}


\begin{figure}
	\includegraphics[width=\textwidth]{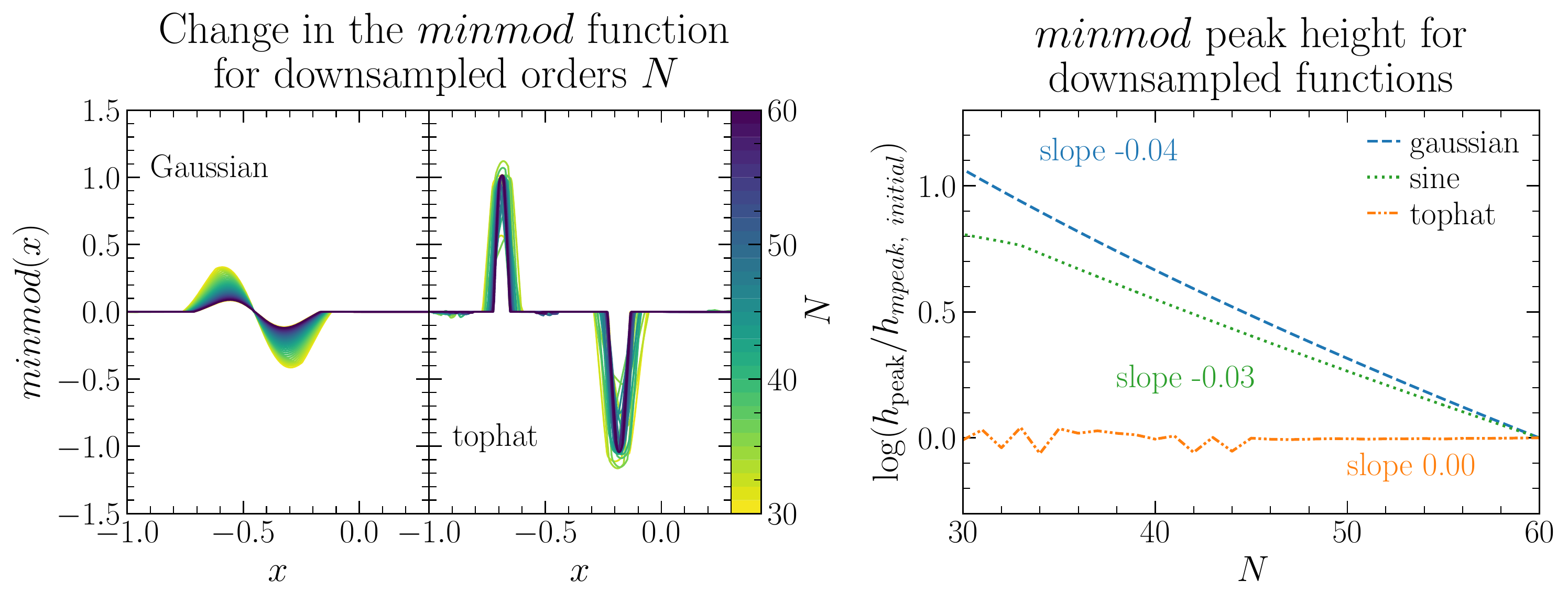}
	\caption{Left panel: behavior of the $minmod$ function under
	a change of the expansion order $N$. The color of each line
	corresponds to an order labelled in the colorbar. For a smooth
	function the $minmod$ peaks converge to 0, while they retain
	constant amplitude for the discontinuous one.
	Right panel: change of minmod peak height with
	the number of modes sampled downward from N=60. 
	Continuous functions show exponential convergence of the
	$minmod$ peak height with linear
	slopes in the log-linear space of ${\rm slope}\approx-0.04$ 
	and $-0.03$. The discontinuous $minmod$ height stays 
	approximately constant with a vanishing slope in the 
	log-linear space.
	}
	\label{fig:solutions:minmods}
\end{figure}

As described in~\cite{Gelb2008} and Section~\ref{sec:Gibbs} 
the $minmod$ function is designed
to approximate the jump function $[f](x)=f(x^+) - f(x^-)$ of a discontinuous function
$f(x)$. Thus, for smooth solutions $minmod$ vanishes for number
of modes $N \rightarrow \infty$. The left panel of 
Fig.~\ref{fig:solutions:minmods} demonstrates this property by 
showing how the $minmod$ function changes as a function of the 
number of modes $N$ for a smooth Gaussian (Eq.~\ref{eq:gaussian},
$x_0 = -0.45$, $\sigma=0.15$) and a discontinuous `tophat' 
function (Eq.~\ref{eq:tophat}, $x_1=-0.7$, $x_2=-0.2$). 
As the number of modes increases, the
Gaussian $mindmod$ peaks decrease in amplitude with the curve 
approaching a flat line. In the `tophat' case, however, the 
$minmod$ peaks do not vanish, retaining their amplitude and 
becoming narrower in shape as the number of modes increases.

The right panel of Fig~\ref{fig:solutions:minmods} further 
quantifies these behaviors by plotting the variation in the
$minmod$ peak height as a function of number of modes $N$ for the
previously introduced tophat, Gaussian and a 
sine $f(x)=\sin(6x)$ in the log-linear space.
The plot clearly illustrates the exponential decay of the peak
height for the smooth functions contrasted with the constant 
amplitude present in the discontinuous case. The slopes resulting
from a linear fit in this space are negative for the smooth functions
and vanishing for the discontinuous one.

In order to create Fig.~\ref{fig:solutions:minmods} we 
initially decompose functions in Eq.~(\ref{eq:gaussian}) and
(\ref{eq:tophat}) with N$=60$ pseudospectral basis of Chebyshev 
polynomials. We then `downsample' the initial representation 
to lower $N=K$,
obtaining spectral expansion coefficients in the new basis 
$a_{k}$ through the following transformation:
\begin{equation}
	a_k = \mathcal{D}_{km} a_m\, ,
\end{equation}
where
\begin{align}
	\mathcal{D}_{km} & = \mathcal{V}_{kj} \mathcal{P}_{jm}\,, \\
	\mathcal{P}_{jm} & = \sum_{m=0}^M T_m(x_j)\, \\
\end{align}
and $[m]=\{0,1\cdots M\}$, $[j]=\{0,1\cdots K\}$ for $K < M$.
The $\mathcal{P}_{jm}$ matrix 
applied to expansion coefficients $a_m$ probes the spectral 
reprojection of order $M$ on the $K+1$ $x_j$ collocation points 
of the new basis of order $K < M$. These collocation-point values
are then translated to spectral coefficients $a_k$ via the inverse 
of the Vandermonde matrix, $\mathcal{V}$. 

The simple test cases presented in Fig.~\ref{fig:solutions:minmods} 
indicate that different smooth functions can potentially exhibit
different decay slopes for the $minmod$ peak heights with increasing
$N$. Thus, in order to find a robust demarcation between purely
smooth and marginally smooth or discontinuous solutions we investigate a set of different analytic functions and compute convergence
slopes of the $minmod$ peak heights for different initial expansion
orders $N$. We then repeat this exercise for a test set of 
discontinuous problems, compare the results and choose a slope
value which would allow us to robustly select smooth cases.

\begin{figure}
	\includegraphics[width=\textwidth]{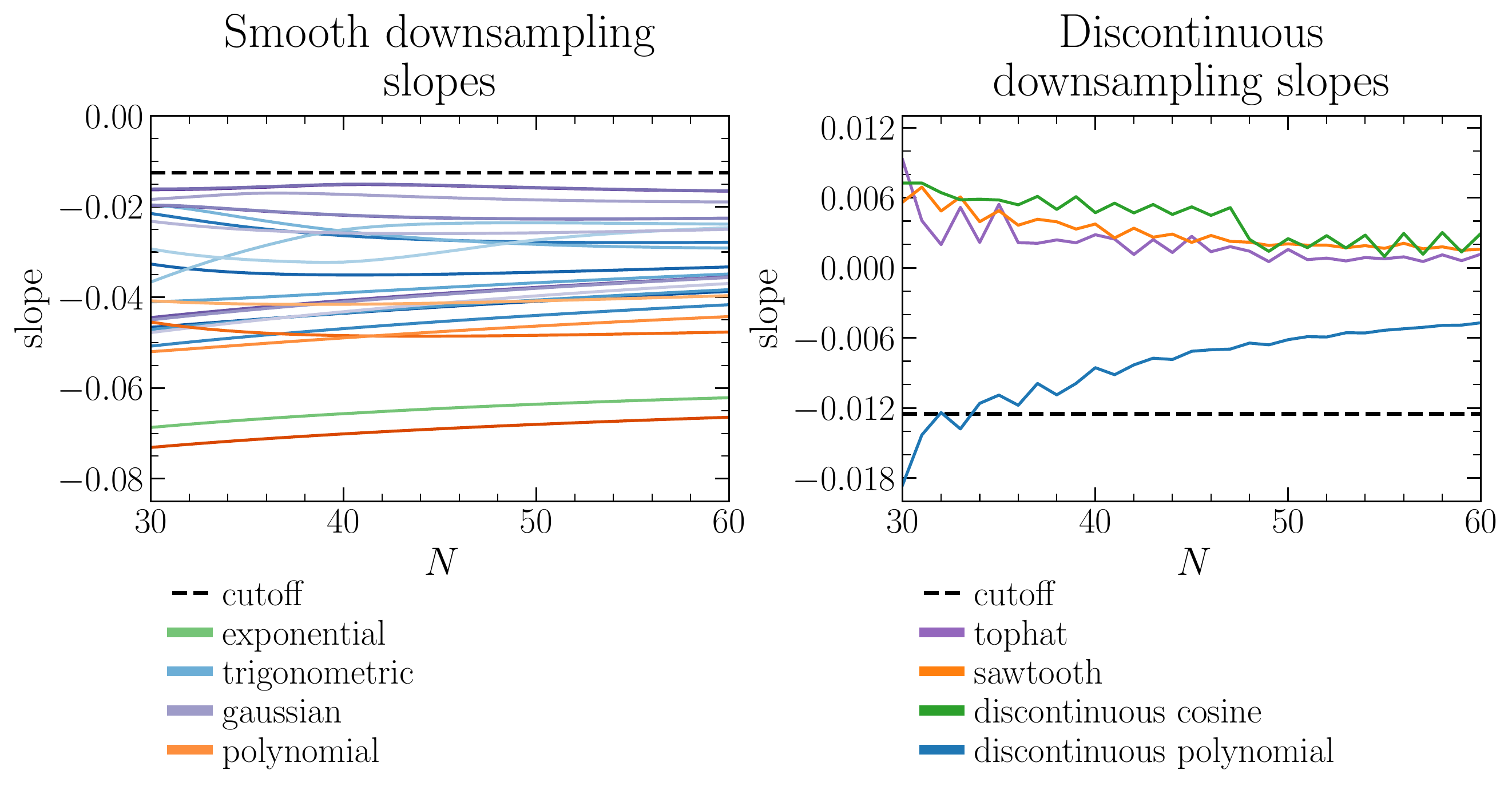}
	\caption{Convergence slopes in the 
	log-linear space for smooth (left panel) and discontinuous
	(right panel) functions explored for different initial
	expansion orders $N$. The investigated range of continuous
	cases is always characterized by a negative slope, while
	majority of discontinuous cases by a positive one. 
	An experimentally suggested demarcation to 
	distinguish purely smooth functions is plotted as 
	a black dashed line.
	}
	\label{fig:solutions:slopes}
\end{figure}

Fig.~\ref{fig:solutions:slopes} presents the results of our 
investigations with the smooth functions tested in the left panel
and the discontinuous ones in the right panel. The set of 25 smooth 
tests contains a variety of exponential, trigonometric, 
Gaussian-like and polynomial functions color-coded by their
respective type labelled in the legend. As presented in the figure
all $minmod$ slopes for the analytic test cases are $<-0.015$ while
majority of the discontinuous ones are $\gtrsim0$. The variation
in slope values of the continuous functions comes from the
variation in the amplitude of the input functions with higher 
amplitudes resulting in increasingly negative slopes.
The polynomial
case is characterized by small magnitudes of discontinuities, which
escape the edge detection prescription and result in different slope
values. After examining our results we choose to use the slope value
of $-0.125$ to robustly label a solution at a given time step as 
smooth.


\vspace{1cm} 
\subsection{Identification and treatment of unresolved steep 
gradients}
\label{sec:smooth-unresolved}


\begin{figure}
	\includegraphics[width=\textwidth]{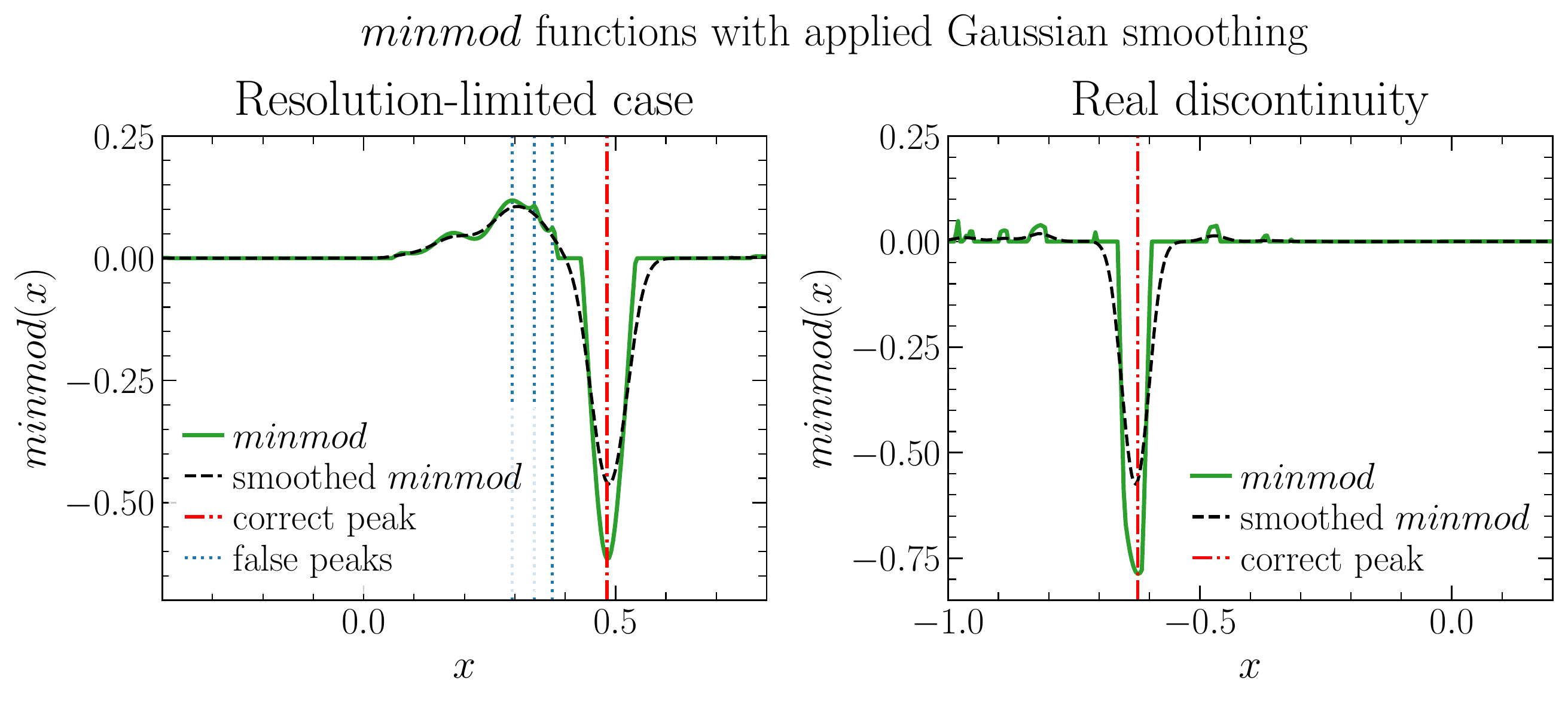}
	\caption{Smoothing procedure which aids in distinguishing
	clear discontinuous cases from resolution-driven 
	oscillations. Left panel: resolution limited case, where the 
	$minmod$ function is polluted by spurious peaks. Initially,
	multiple peaks are recognised, which are then rejected
	in secondary search after Gaussian smoothing is applied.
	Right panel: true discontinuity case where no spurious
	peaks are identified in which case the number of detected
	peaks does not change after applying gaussian smoothing.
	}
	\label{fig:solutions-2}
\end{figure}

When the solution fails the smoothness test at a given
time step we then search for the potential shock locations by
identifying local extrema in the $minmod$ function. As described
earlier in Section~\ref{sec:detection-challenges} the unresolved
gradients create spurious extrema revealed by the automated
detection algorithm. Fig.~\ref{fig:solutions-2} illustrates
this issue by plotting a comparison between
the $minmod$ functions resulting from a resolution limit and 
a genuine discontinuity. The former is shown in the left panel,
where three spurious peaks in the function are visible, marked in blue 
dotted lines. Location of the steep gradient is marked in
red dash-dotted vertical line. This case is in visible contrast
with the right panel, where the $minmod$ function clearly exhibits
a single, distinct peak and no spurious instances are identified. 

In order to remove unwanted extrema we developed a three-step
heuristic procedure performed as an additional part of the edge 
detection routine. It begins with an initial search for extrema
as the first step, which identifies all potential peaks, including
those marked with dotted lines in the left panel of 
Fig.~\ref{fig:solutions-2}. 

In the second step, each of the initially 
found extrema are treated separately to verify their identification.
For a given suggested jump location we construct a Gaussian smoothing
kernel of width equal to the distance between the two collocation 
points $x_i$ and $x_{i+1}\,$, neighboring the extremum: 
\begin{equation}
	k(x-x^\prime) = \exp(\frac{-(x-x^\prime)^2}{2\omega^2})\, , 
	\label{eq:smoothing-kernel}
\end{equation}
where $\omega = (x_{i+1} - x_{i})/2$. 

We then convolve the $minmod$ with this kernel: 
\begin{equation}
	minmod(x)_{\rm smoothed}= 
	\int^{1}_{-1}minmod(x^\prime)k(x-x^\prime)
	dx^\prime
\end{equation}
to wash away any small peaks resulting from the 
inability of the fixed grid of collocation points to capture steep
gradients in smooth functions.

As a final step we repeat the search for extrema
restricting their width to no more than twice the distance between
the adjacent collocation points. If a given extremum is 
identified in the repeated search it is accepted as a genuine 
feature, otherwise labelled as a spurious one. As illustrated in 
Fig.~\ref{fig:solutions-2}, the narrow features removed
by smoothing do not show up in the search while the wide extremum 
is discarded by the maximum width criterion. Thus, only the 
extremum corresponding to the location of the
steepest gradient in the solution is retained in this example.

In the edge detection routine, once the artefacts are recognized, 
the solution at a given time step is labelled as resolution-limited case.
Such cases also require mollification treatment in order to 
remove Gibbs oscillations and we choose to use a continuous
optimal mollifier as introduced in Eq.~(4.2) of \cite{Tanner2006}. 
This operation
allows us to recover the solution without Gibbs oscillations, 
simultaneously preserving the steep gradient without excessive
smoothing.

If, on the other hand no spurious extrema are confirmed in 
the repeated search, the case is labelled as genuinely 
discontinuous and the edge location found in the first step
is correct (e.g. right panel of Fig.~\ref{fig:solutions-2}).
Such a solution is then mollified using a prescription
for one-sided mollifiers from Eq.~(28) in \cite{Piotrowska2019} to
recover the Gibbs-free solution along with its sharp discontinuity
(shock) feature. 

\vspace{0.8cm}
\section{Application: 1D inviscid Burgers' equation}
\label{sec:burgers}

We apply our improved shock detection techniques to a toy problem
involving dynamical development of discontinuities from smooth
solutions. We solve the 1D inviscid Burgers' equation from 
Eq.~(\ref{eq:inviscid-burgers}) with a $N=60$ pseudospectral 
Chebyshev basis under periodic boundary conditions.
As our initial conditions we choose a unit height Gaussian
function in Eq.~(\ref{eq:gaussian}) with $x_0=0$ and $\sigma=0.15$. 
To evolve the nonlinear system
in time we use the filter defined in Eq.~(75) of \cite{Miller2016} 
applied to the time derivative to ensure stability of the numerical
scheme. We allow the system to evolve until $t=3.0$ and then treat
the Gibbs-affected solution with the edge detection and mollification
techniques as described earlier in Sec.~\ref{sec:our-improvements}.

\begin{figure}
	\includegraphics[width=\textwidth]{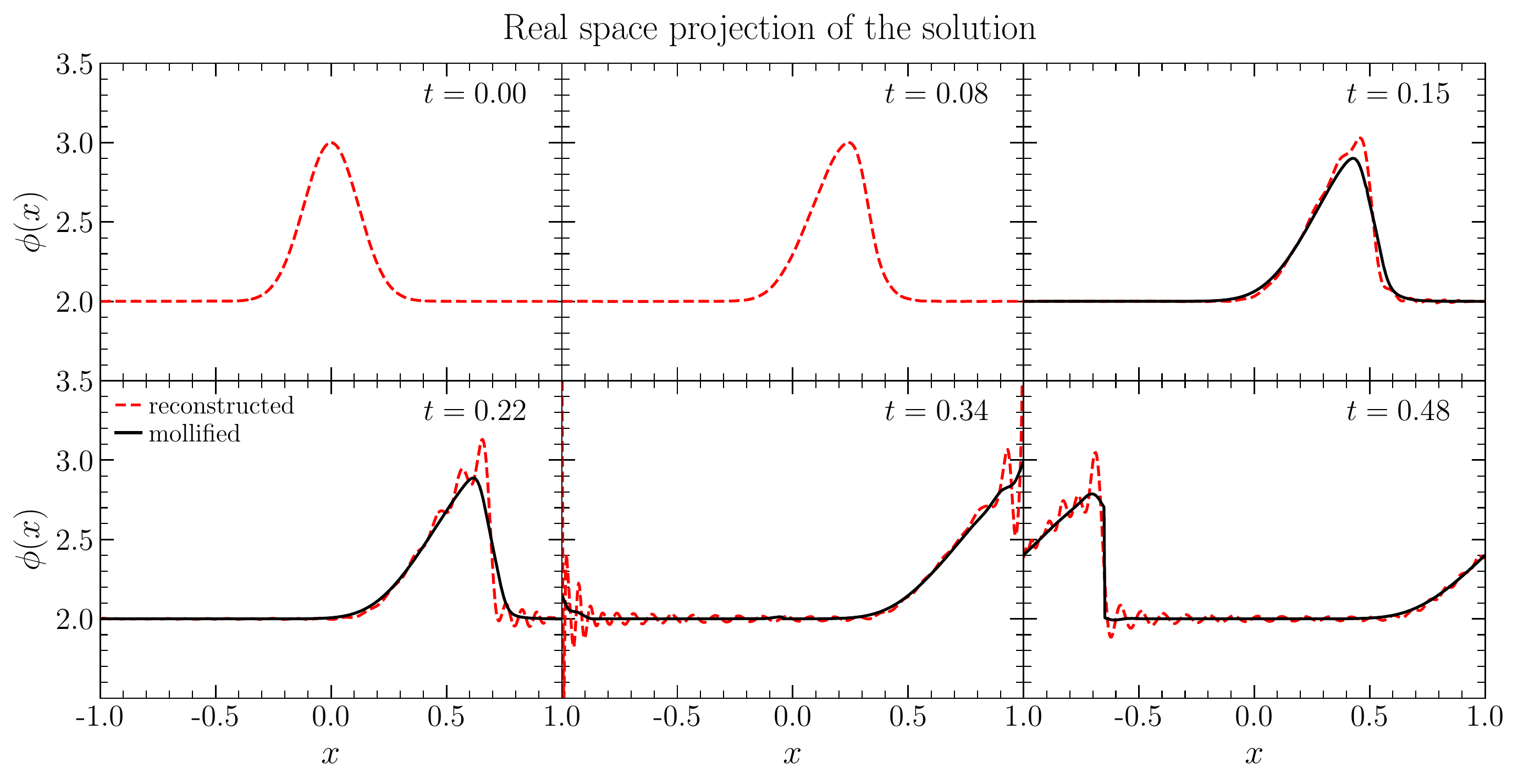}
	\caption{Snapshots from pseudospectral solution of the
	1D inviscid Burgers' equation post processed with our
	improved shock detection and appropriate choice of 
	mollification technique. The red dashed lines show
	the raw solution while the black solid ones our 
	post-processed result. The time progression is left to right
	and top to bottom.
	The first two snapshots do not feature mollification because
	the solution was recognised as smooth. The solid curves show
	a solution free from spurious oscillations, smoothly 
	evolving to form a discontinuity captured by our method. 
	}
	\label{fig:burgers-evolution}
\end{figure}

Fig.~\ref{fig:burgers-evolution} presents snapshots of our solution
at different times with the spectral reconstruction plotted in 
red dashed line and the mollified reconstruction overlaid as a black
solid curve. The time progression is from left to right and from
top to bottom with the initial conditions shown in the upper left
corner and a $t=0.48$ snapshot in the right bottom one. The black 
curve is missing in the first two snapshots because the solution
was recognized as smooth within the resolution limit and didn't
require post-processing to remove Gibbs-like oscillations. The 
third snapshot demonstrates how the resolution affects the 
vicinity of steep gradients and how well continuous mollifies
cope with removing spurious information introduced locally into
the solution. Our improved edge detection techniques are also capable
of correctly identifying steep gradients located across the domain
boundaries (middle bottom panel) allowing for robust treatment of
marginally continuous cases located at the least favorable location
within the computational domain. The last panel shows a solution
with a~clearly developed discontinuity where the Gibbs oscillations
pollute the entire domain. The edge detection identified a 
single extremum in the $minmod$ function and hence one-sided 
mollification was successfully applied to recover the Gibbs-free 
solution while preserving its discontinuous character.

\vspace{0.8cm}
\section{Summary and conclusions}
\label{sec:conclusions}

Pseudospectral methods are a powerful means of solving PDEs which
enjoy global exponential convergence for smooth problems.
When applied to treat discontinuous functions, however, they lose 
their main advantages and suffer spurious oscillations in the 
solution which deems them unsuitable for astrophysical applications
involving for e.g. fluid shocks. Over the years, multiple authors
developed post-processing techniques allowing for a successful
removal of the Gibbs phenomenon which we have previously tested
in practice and applied to a toy advection problem.

In this work we focused on discontinuities which develop dynamically
from smooth solutions, highlighting the challenges associated with
their correct detection and treatment. We showed how steepening 
gradients lead to Gibbs-like fluctuations concentrated close to the
discontinuities and how they influence the edge detection procedure
causing problems in automated identification of the shock creation.

We then propose a heuristic solution to those problems, which
allow for a robust automated identification of solutions 
requiring post-processing. This step is an inseparable part
of any (pseudo)spectral 
shock capturing scheme, regardless of the choice of
oscillation removal technique. Thus, in this work we offer
solutions which can potentially benefit a wider community,
encouraging further application of pseudospectral schemes
to solving discontinuous problems.

Our improvements to jump identification 
utilize the rate of convergence of the $mimod$ function to 
differentiate perfectly smooth cases from the marginally smooth ones.
We also suggest the use of Gaussian smoothing kernels on the
$minmod$ function to identify solutions affected by the resolution 
limit. For these cases we then suggest treatment with symmetric 
adaptive mollifiers, while for the genuine shocks we recommend
the one-sided kind. We tested this prescription on a range of
smooth functions involving steep gradients located at random
positions within the domain. The procedure was very successful
at tackling steep gradients and genuine discontinuities
across the whole domain, and suffered less problems close to the 
domain boundaries. 

The combination of Gaussian smoothing 
and repeated peak search allows us to avoid problems
akin to those highlighted in \cite{Piotrowska2019} such as
spurious, narrow peaks originating close to the domain
boundaries.
Nevertheless, the method still suffered from occasional
misidentification of genuine discontinuities as steep
gradients, especially in cases when the discontinuity 
fell exactly at the domain boundary. Aware of this numerical
weakness, we recommend exercising caution when dealing
with discontinuities located close to the domain boundaries.

Finally, we solved a 1D inviscid Burgers' equation to test the
performance of our implemented methods in more realistic settings. Our
framework was able to automatically determine when the gradients in
the solution become steep enough to require post processing and treat
the genuine discontinuity once it develops. Although we only tackled
the inviscid Burger's equation with a single set of initial
conditions, we argue that our approach extends to arbitrary
systems. In particular, one must introduce length and height scales,
which need to be incorporated into the convergence rates shown in
figure \ref{fig:solutions:slopes}.

Encouraged by our progress with dynamically developing
discontinuities for test cases we will apply 
our pseudospectral Chebyshev framework to solving 
real shock problems in the future. Benefiting from 
the apparent robustness of the scheme we would also
like to extend it to two dimensions in order to 
tackle astrophysical problems with accuracy higher
than the popular finite volume methods.

\section*{Acknowledgements}

We would like to thank Erik Schnetter, Josh Dolence and Chris Fryer
for their support and constructive comments. JMP would also like to
thank Roberto Maiolino and Asa Bluck for their excellent mentorship
and valuable suggestions which increased the quality of this paper.

We gratefully acknowledge support from the U.S. Department of Energy
(DOE) Office of Science and the Office of Advanced Scientific
Computing Research via the SciDAC4 program and Grant DE-SC0018297,
from the U.S. NSF grant AST-174267, and from the U.S. DOE through Los
Alamos National Laboratory (LANL). This work used resources provided
by the LANL Institutional Computing Program. Additional funding was
provided by the LDRD Program and the Center for Nonlinear Studies at
LANL under project number 20170508DR. LANL is operated by Triad
National Security, LLC, for the National Nuclear Security
Administration of the U.S. DOE (Contract
No. 89233218CNA000001).

This article is cleared for unlimited release, LA-UR-19-28857.

We are grateful to the countless developers contributing to open
source projects on which we relied in this work, including Python
\citep{rossumPythonWhitePaper}, numpy and scipy
\citep{numpy,scipyLib}, and Matplotlib \citep{hunterMatplotlib}.





\bibliographystyle{elsarticle-num}

\bibliography{paper-draft}

\end{document}